\documentclass[12pt]{article}
\usepackage{amsfonts,amsmath,bbm}
\def\proof{\noindent{\bf Proof:}\hskip10pt}        
\def\QED{\hfill $\Box$}

\font\tenmath=msbm10 scaled 1200
\font\sevenmath=msbm7 scaled 1200
\font\fivemath=msbm5 scaled 1200
\newfam\mathfam \textfont\mathfam=\tenmath
\scriptfont\mathfam=\sevenmath \scriptscriptfont\mathfam=\fivemath

\begin{document}
\def \\ { \cr }
\def \R{\mathbb{R}}
\def\N{\mathbb{N}}
\def\E{\mathbb{E}}
\def\P{\mathbb{P}}
\def\Z{\mathbb{Z}}
\def\D{\mathbb{D}}
\def\C{\mathbb{C}}
\def\da{^{\downarrow}}
\def \e{{\rm e}}
\def \p{{\cal P}}
\def \s{{\cal S}}
\def \g{{\cal G}}
\newtheorem{theorem}{Theorem}
\newtheorem{definition}{Definition}
\newtheorem{proposition}{Proposition}
\newtheorem{lemma}{Lemma}
\newtheorem{corollary}{Corollary}
\centerline{\LARGE \bf Different Aspects of a Model for }
\vskip 2mm
\centerline{\LARGE \bf Random Fragmentation Processes}

\vskip 1cm
\centerline{\Large \bf Jean Bertoin}
\vskip 1cm
\noindent
\centerline{\sl Laboratoire de Probabilit\'es et Mod\`eles Al\'eatoires}
\centerline{\sl  et Institut
universitaire de France}
\centerline{\sl Universit\'e Pierre et Marie Curie, et C.N.R.S. UMR 7599}
\centerline{\sl 175, rue du Chevaleret} 
\centerline{\sl F-75013 Paris, France}
\vskip 15mm

\noindent{\bf Summary. }{\small This text surveys different
probabilistic aspects of a model which is used to describe the evolution of an
object that falls apart randomly as time passes. Each point of view
yields useful techniques to establish properties of such random fragmentation
processes.}
 \vskip 5mm
\noindent
{\bf A.M.S. Classification.}  {\tt 60J80, 60G18}
\vskip 3mm
\noindent{\bf e-mail.} {\tt jbe@ccr.jussieu.fr}

\section{Introduction}
Fragmentation is natural phenomenon that can be
observed at a great variety of scales. To give just
a few examples, we may think of stellar fragments and meteoroids in Astrophysics,
fractures and earthquakes in Geophysics, crushing in the mining industry,  breaking of
crystals in Crystallography, degradation of large polymer chains in Chemistry, fission of
atoms in Nuclear Physics, fragmentation of a hard drive or files in Computer Science, ...
In this text, we will be interested in situations where this phenomenon
occurs randomly and repeatedly as time
passes. Typically, we may imagine the evolution of blocks of mineral in
a crusher. 

Over many years, the importance of the activity  on fragmentation
in Physics (see e.g. \cite{Ald, BCP} and the references therein) has not been
reflected in Probability Theory where the interest has been more irregular. The first
significant works concerning its  probabilistic aspects (strong laws) are due to Kolmogorov
\cite{Ko} and his student Filippov \cite{Fi}.  The systematic study of general
fragmentation processes is quite recent, although some well-established areas of 
Probability Theory, such as branching random
walks and multiplicative cascades, are clearly relevant to investigate a large family
of fragmentation processes.

This text is intended as a survey on a stochastic model for fragmentation processes which
is characterized by a few parameters. We shall present various point of views on this
model, from which different properties depending
on the values of the parameters are derived. In order to deal with models that can be
studied mathematically, we are led to make hypotheses that may look at
first sight somewhat stringent,
but which are often fulfilled in applications. First, we suppose that the system has a
memoryless evolution, i.e. its future only depends on its present state and not on its
past. In other words, the system enjoys the Markov property. In particular, this excludes
the possibility that an object might be more fragile (i.e. more likely to split) due
to former shocks. Second, we assume that each fragment can be characterized by a real number
that can be thought of as its size. This impedes to consider geometrical properties like
the shape of a fragment. Third, we shall suppose that the branching property is fulfilled,
in the sense that fragments split independently one of the other, or in other words, that
the evolution of a given fragment does not depend on its environment. Finally, we shall
assume that the process enjoys self-similarity, that is that the law of
a fragmentation process started from a unique fragment of size $r>0$
can be reduced by proper rescaling in space and time, to that when $r=1$.

In order to give a formal definition, we introduce the state space of
decreasing numerical sequences which tend to $0$,
$$\s\,:=\,\left\{{\bf s}=(s_1,s_2,\ldots): s_1\geq s_2\geq\ldots \geq0
\hbox{ and }\lim_{n\to\infty}s_n=0\right\}\,.$$
A generic configuration ${\bf s}$  
should be thought of as the ranked sequence of the sizes of some object
that has been split.
The space $\s$  is
endowed with the uniform distance, so $\s$ is a complete metric separable space.

Let $X=(X(t), t\geq0)$ denote a random process 
such that 
$$X(t)=(X_1(t),X_2(t), \ldots)\in \s\qquad \forall t\geq0,\ a.s.$$
We suppose that $X$ is continuous in probability and Markovian. For every
$r\geq0$, we write 
$\P_r$ for its law started from the configuration $(r,0,\ldots)$,
and we further assume that
for every $t\geq0$, $\P_1(X_1(t)\leq 1)=1$ and
$\P_0(X_1(t)=0)=1$. In other words, the sizes of the pieces resulting at time $t$ from an
object with initial size $1$ cannot exceed $1$, and an object with $0$ size cannot produce
pieces with positive size.

\begin{definition} \label{D1}  
We call $X$ a {\rm self-similar fragmentation} if :

\noindent {\rm (i)} There exists $\alpha\in \R$, called the {\rm index of
self-similarity}, such that for every $r>0$, the distribution under $\P_1$ of the
rescaled process
$\left(rX(r^{\alpha}t), t\geq0\right)$ is $\P_r$. 

\noindent {\rm (ii)} For every ${\bf s}=(s_1,s_2,\ldots)\in\s$,
if $\left(X^{(i)},i\in \N\right)$ is a sequence of independent processes
such that $X^{(i)}$ has the law $\P_{s_i}$, and if 
$\tilde X(t)$ denotes the decreasing rearrangement of the family
$\left(X^{(i)}_n(t): i,n \in \N\right)$, then $\tilde X
=({\tilde X}(t), t\geq0)$ is a version of $X$ started from the
configuration ${\bf s}$.

\end{definition}

The self-similarity and branching assumptions (i-ii) enable us to focus on the
special case when the fragmentation starts
from a single fragment with unit size, without losing generality.
Thus, in the sequel, we whall implicitly work under the probability
measure $\P:=\P_1$.

It is intuitively obvious that the behavior of a self-similar
fragmentation should depend crucially on the sign of the index of
self-similarity. Roughly, since fragments get smaller as time passes,
 the rate of dislocations decreases  when the index
is positive, whereas it increases when the index is negative. 
In the critical case when $\alpha=0$, the fragmentation process is
called {\it homogeneous}, as the rates at which fragments split do not depend
of their sizes. Homogeneous fragmentations are of course easier to study,
and more results are known for homogeneous fragmentations than for
general ones.

In the case when splits occur discretely, that is when each fragment
remains unchanged for some random time and then splits, fragmentation processes
can be described in terms of
branching Markov chains in continuous time. Further, they can be endowed
with a discrete genealogical structure, which enables us to develop their
study using general techniques from the theory of multiplicative
cascades.  Plainly, this approach fails when splits occur
continuously,  that is when each fragment breaks down instantaneously.
In the sixth section of this paper, we shall circumvent this fundamental difficulty by
performing a spatial discretization. Specifically, Kingman's theory of 
exchangeable random partitions provides the right framework for studying
a fundamental special class of fragmentation processes, 
called exchangeable. This approach enables us
to reveal the fine structure of exchangeable fragmentation
and to extend results which were proven in the discrete setting to this general
setting. Finally, the last section of this paper is devoted to a rather informal
discussion of the duality relation involving time-reversal, which exists between certain pairs of fragmentation and coalescent processes.

\section{Fragmentations as branching Markov chains}
\subsection{Construction of fragmentation chains}
A self-similar fragmentation process is called a (self-similar fragmentation) {\it chain}
if the first dislocation time
$$T:=\inf\left\{t>0: X(t)\not=(1,0,\ldots)\right\}$$
is strictly positive a.s. (recall that we implicitly assume 
that at the initial time there is a
single fragment with unit size).
In this situation, the Markov property forces
$T$ to have an exponential distribution, say with parameter $c\geq 0$. Excluding
implicitly the degenerate case when $c=0$ (i.e. $T=\infty$ a.s.), 
 we may and will henceforth focus on the case when
$c=1$ for the sake of simplicity, as the general case can be reduced to the former by a
linear time-substitution. We write $\nu$ for the distribution of $X(T)$ under $\P_1$.
So $\nu$ is a probability measure on $\s$ with $\nu(\{(1,0,\ldots)\})=0$, which we
call the {\it dislocation law} of $X$.

We can then think of $X$ as the evolution of a non-interacting particle
system in $]0,\infty[$, in which each fragment with positive size
is viewed as a particle and different particles have independent evolutions.
Specifically,  a particle with size $r>0$ lives for an
exponential time with parameter $r^{\alpha}$. Then it disappears
and is replaced by a cloud of smaller particles, say $(r_1,\ldots)\in\s$,
such that the sequence of ratios $(r_1/r, \ldots)$ has the fixed law
$\nu$ and is independent of the lifetime of the particle
$r$. 
This description makes sense only
when $r>0$; however by assumption, the possible children
of a particle with size $0$ have all size $0$. Particles with size $0$ play no role,
and the evolution is thus well-defined in all cases.
It should be intuitively clear from this description
that the law of the system is entirely determined by the dislocation
measure $\nu$ and the index of self-similarity $\alpha$.

Conversely, given  a real number $\alpha\in\R$ and a probability measure $\nu$ on $\s$
which fulfils mild conditions,
it is easy to construct a self-similar fragmentation chain
with index $\alpha$ and dislocation law $\nu$.
Roughly, the idea is to focus first on particles with size at least $\varepsilon>0$,
and to built the restricted system using the standard theory of Markov chains in
continuous time. Then one observes a compatibility property for
different values of the threshold  parameter $\varepsilon>0$, which allows us to take
a projective limit as $\varepsilon\to0$.

More precisely, assume that
$\nu(\{{\bf s}\in\s: s_1>1\})=0$,  $\nu(\{1,0,\ldots)\})=0$, and
\begin{equation}\label{eqSF1}
\int_{\s}\#\left\{i: s_i>\varepsilon\right\}\nu(d{\bf
s})\,<\,\infty\qquad \hbox{for every } \varepsilon>0.
\end{equation}
For every $\varepsilon>0$, one can then construct
\footnote{Condition (\ref{eqSF1}) is needed to prevent a
possible explosion, as otherwise it could happen
that the dynamics would create an infinite number of atoms of size $>\varepsilon$
in finite time.}
a branching Markov chain in continuous time, say $X^{(\varepsilon)}$, with
values in the space of finite atomic measures on $]\varepsilon,\infty[$
and governed by the following transitions. When the process starts from
say $m=\sum_{i=1}^{n}\delta_{r_i}$ with $r_1\geq r_2\geq \ldots \geq
r_n>\varepsilon$, the first jump occurs at an exponential time
with parameter $\sum_{i=1}^{n}r_i^{\alpha}$. The configuration
immediately after the jump is independent of the exponential time
and has the same distribution as the point measure obtained from
$m$ by picking an atom $r_i$ at random with probability proportional
to $r_i^{\alpha}$, and replacing it by a random family of atoms,
say $(r_i s_1,\ldots, r_i s_k)$
where ${\bf s}=(s_1,\ldots)$ has the law $\nu$
and $k=\max\{j\in \N: r_is_j>\varepsilon\}$.
It is then immediately checked from calculations of branching rates that
for every
$\varepsilon>\eta>0$, the restriction of the point measure process
$X^{(\eta)}$ to $]\varepsilon,\infty[$ is a version of
$X^{(\varepsilon)}$. By Kolmogorov's extension theorem, this ensures the
existence of a unique law of a process $X=(X(t), t\geq0)$ with values in
the space of Radon point measures on $]0,\infty[$, such that the
restriction of $X$ to $]\varepsilon,\infty[$ is a version of
$X^{(\varepsilon)}$ for every $\varepsilon>0$. The natural
identification of a Radon point measure with the ranked sequence of its
atoms enables us to view $X$ as a self-similar fragmentation chain
with the desired dynamics.

We stress that in general, despite of the terminology, $X$ is {\it not} a continuous time
Markov chain. Indeed, even though the lifetime of a fragment is exponentially
distributed with a positive parameter, the process
may create infinitely many fragments in a finite time (but of course, only
finitely many of them have a size greater than $\varepsilon$ for any
fixed $\varepsilon>0$), and the infimum of the lifetimes of an infinite
family of fragments can be zero.  

In the special case $\alpha=0$, the lifetime of each 
fragment is a standard exponential variable, independently of the size of the fragment,
and  we say that the fragmentation chain is {\it homogeneous}. In this situation, there is
a natural connection with branching random walks in continuous times (cf. Uchiyama
\cite{U}). Specifically, consider a homogeneous fragmentation chain $X$
with a dislocation law $\nu$ that
charges only the sub-space of sequences ${\bf s}=(s_1, \ldots)$ with $s_k=0$ for $k$
sufficiently large. Then the process
$$Z^{(t)}(dx):=\sum \delta_{-\ln X_i(t)}(dx)\,,$$
where the sum is taken over the fragments with strictly positive size,
is a branching random walk. More precisely, its branching measure is
the image of $\nu$ by the map $x\to-\ln x$.
This elementary connection  has a number of interesting consequences as it essentially
reduces the study of the class of homogeneous fragmentations associated to a
dislocation laws charging only finite configurations, to that of branching random walks on
$]0,\infty[$,  for which a great deal of results are known.

The infinitesimal generator ${\bf G}$ of a self-similar fragmentation
chain has a simple and useful expression for so-called additive and
multiplicative functionals, which is easily obtained from the construction
described above in terms of the jump rates.

\begin{proposition}\label{P1} 
\noindent{\rm (i)} Consider a measurable function $f:[0,\infty[\to\R$ with
$f=0$ on some neighborhood of $0$, and define an additive functional
$A:\s\to\R$ by
$$A({\bf s})\,=\,\sum_{i=1}^{\infty}f(s_i)\,, \qquad {\bf s}=(s_1, s_2, \ldots)\,.$$
Then for every ${\bf x}=(x_1, \ldots)\in \s$, we have
$${\bf G}A({\bf x})\,=\,\sum_{i=1}^{\infty}{\bf
1}_{\{x_i>0\}} x_i^{\alpha}\int_{\s}(A(x_i{\bf s})-f(x_i))\nu(d{\bf s})\,.
$$

{\noindent \rm (ii)} Consider a measurable function
$g:[0,\infty[\to]0,\infty[$ with $g=1$ on some neighborhood of $0$, and
define a multiplicative functional $M:\s\to]0,\infty[$ by
$$M({\bf s})\,=\,\prod_{i=1}^{\infty}g(s_i)\,, \qquad {\bf s}=(s_1, s_2, \ldots)\,.$$
Then for every ${\bf x}=(x_1, \ldots)\in \s$, we have
$${\bf G}M({\bf x})\,=\,\sum_{i=1}^{\infty}{\bf 1}_{\{x_i>0\}} x_i^{\alpha}{M({\bf x})\over
g(x_i)}\int_{\s}(M(x_i{\bf s})-g(x_i))\nu(d{\bf s})\,.$$
\end{proposition}

\subsection{Some analytic expressions and
 Malthusian hypotheses}

Throughout the rest of this section, we will be working
with  a self-similar fragmentation
chain $X$ with index $\alpha\in\R$ and dislocation law $\nu$,
where $\nu$ is a probability measure of $\s$ satisfying the
conditions of the preceding section. To avoid uninteresting discussions,
we shall not consider the situation when $s_i=0$ or $1$ for all $i\in\N$,
$\nu(d{\bf s})$-a.s., since the fragmentation chain can then be identified as a branching
process in continuous time.

  For future purposes, it is convenient to introduce
 the notation
$$\underline p:=\inf\left\{p>0:  \int_{\s}\sum_{i=1}^{\infty}s_i^p\nu(d{\bf s})<\infty\right\}$$
(with the convention $\inf \emptyset = \infty$), and to assume from now
on that $\underline p<\infty.$
Then we define for every
$p\geq\underline p$
\begin{equation}\label{eq2}
\kappa(p):=\int_{\s}\left(1-\sum_{i=1}^{\infty}s_i^p\right)\nu(d{\bf s})\,.
\end{equation}
Note that  our assumptions ensure that $\kappa$ is always a continuous strictly increasing
function on
$]\underline p,\infty[$;
$\kappa(\underline p-)$ may be finite or equal to $-\infty$.

We then make the fundamental :

\vskip 2mm
\noindent
{\bf Malthusian Hypotheses.} {\it There exists a (unique) solution $p^*> \underline p$
to the equation
$\kappa(p^*)=0$, which is called the  {\rm  Malthusian exponent}. Furthermore the
integral
$$\int_{\s}\left(\sum_{i=1}^{\infty}s_i^{p^*}\right)^p\nu(d{\bf s})$$ 
is finite for some $p>1$.}

\vskip 2mm
\noindent We stress that the Malthusian hypotheses are quite weak. For instance, if 
$s_1<1$, $\nu(d{\bf s})$-a.s. (which means that the fragments resulting from the
dislocation of a particle are always strictly smaller than the initial particle), then
we have by dominated convergence that $\lim_{p\to\infty}\kappa(p)=1$.
If moreover we can find some $p>\underline p$ such that $\kappa(p)<0$
(this occurs for instance whenever $\int_{\s}\#({\bf s})\nu(d{\bf s})\in ]1,\infty[$
where $\#({\bf s}):={\tt Card}\{i: s_i>0\}$, then the fact that $\kappa$ is s a
continuous and strictly increasing ensures the existence of the Malthusian
parameter.

Throughout the rest of this work, the Malthusian hypotheses
will always be taken for granted. The function $\kappa$ 
and the Malthusian exponent $p^*$  play a crucial role in the study of the asymptotic
behavior of fragmentation chains, as we shall see soon.

\section{Fragmentation chains and multiplicative cascades}
In this section, we point at a representation of 
a fragmentation chain as an infinite  tree with
random  marks. This representation
can be viewed as a different parametrization of the process,
in which the natural time is replaced by the generation of the different particles.

\subsection{Genealogical coding of fragmentation chains}
We start by introducing some notation. We consider the infinite tree
$${\cal U}\,:=\,\bigcup_{n=0}^{\infty}\N^n\,,$$ with the convention
$\N^0=\{\emptyset\}$. 
In the sequel ${\cal U}$ will often be referred to as the
{\it genealogical tree}; its elements are called {\it nodes} (or sometimes also individuals) and the
distinguished node $\emptyset$ the {\it root}. Nodes will be used to label the particles produced by a
fragmentation chain.
 For each $u=(u_1,\ldots,u_n)\in{\cal U}$, we
call {\it generation} of $u$ and write $|u|=n$, with the obvious convention
$|\emptyset|=0$. When $n\geq 1$ and
$u=(u_1,\ldots,u_n)$, we call
$u-:=(u_1,\ldots,u_{n-1})$ the parent of $u$. Similarly, for every $i\in \N$ we write
$ui=(u_1,\ldots,u_n,i)\in\N^{n+1}$ for the $i$-th child of $u$.
Finally, we call {\it mark} any map from ${\cal U}$ to some (measurable) set.

Now, we associate to each trajectory of the fragmentation chain a
mark on the infinite tree
${\cal U}$. The mark at a node $u$ is the triple $(\xi_u, a_u, \zeta_u)$ where $\xi_u$ is
the size, $a_u$ the birth-time and $\zeta_u$ the lifetime of the particle with label $u$.
More precisely, the initial particle
corresponds to the ancestor
$\emptyset$ of the tree ${\cal U}$, and the mark at $\emptyset$ is the
triple $(1,0,\zeta_{\emptyset})$ where
$\zeta_{\emptyset}$ is the lifetime of the initial particle; in particular
$\zeta_{\emptyset}$ has the standard exponential law.
The nodes at the first generation are used as the labels
of the particles arising from the first split. Again, the mark associated to each
 node
$i\in \N^1$ at the first generation, is the triple $(\xi_i, a_i,\zeta_{i})$, where
$\xi_i$ is the size of the $i$-th child of the ancestor,
$a_i=a_{\emptyset} +\zeta_{\emptyset}$ (the birth-time of a child particle coincides with
the death-time of the parent), and
$\zeta_{i}$ stands for the lifetime of the $i$-th child. And we iterate the same 
construction with each particle at each generation.

Clearly, the description of the dynamics of fragmentation
entails that its genealogical coding
also enjoys the branching property.
Specifically, we have the following recursive description : 

\begin{proposition}\label{PSF2'} There exists two independent families of i.i.d. variables
indexed by the nodes of the genealogical tree,
$\left(\tilde \xi_{u\bullet}, u\in{\cal U}\right)$ and $\left({\bf e}_{u,\bullet}, u\in
{\cal U}\right)$, where each $\tilde \xi_{u\bullet}:=(\tilde
\xi_{u_1,\ldots,u_n,i},i\in\N)$ is distributed according to the law $\nu$, and each
${\bf e}_{u,\bullet} =({\bf e}_{ui}, i\in \N)$ is a sequence of i.i.d.
standard exponential variables, and such that the
following holds:

\noindent
Given the marks $\left((\xi_v,a_v,\zeta_v), |v|\leq n\right)$ of the
first
$n$ generations, the marks at generation $n+1$ can be expressed in the
form
$$(\xi_{ui},a_{ui},\zeta_{ui})
\,=\,(\tilde \xi_{ui}\xi_{u},a_{u}+\zeta_{u},\xi_{ui}^{-\alpha}{\bf e}_{ui}
)\,,$$ 
where $u=(u_1,\ldots,u_n)$ and $ui=(u_1,\ldots,u_n,i)$ is the $i$-th child of $u$.
\end{proposition}

Proposition \ref{PSF2'} shows that the sizes at nodes $(\xi_u, u\in{\cal U})$ define a
so-called {\it multiplicative cascade}; see the pioneer works of Mandelbrot \cite{Mdb},
Kahane and Peyriere
\cite{KP},  Mauldin and Williams \cite{MW}; see also Liu
\cite{Liu1} for further references. Although this multiplicative cascade alone does not
enable us to recover the fragmentation chain (we also need the information on
birth-times and lifetimes as it is amplified in Proposition \ref{PSF2"} below), classical
notions and  results in this field have a key role in the study of fragmentation chains.
It should be intuitively obvious that one can express the fragmentation chain at time $t$
in terms of the particles  which are alive at time $t$, i.e. which are born at or before
$t$ and die after $t$.

\begin{proposition}\label{PSF2"}
With probability one, for every $t\geq0$ and every measurable function
$f:[0,\infty[\to[0,\infty[$ with $f(0)=0$, there is the
identity
$$\sum_{i=1}^{\infty}f(X_i(t))\,=\,\sum_{u\in{\cal U}}{\bf 1}_{\{a_u\leq t < a_u+\zeta_u\}}
f(\xi_u)\,.$$
\end{proposition}

\proof We have to check that all the fragments with positive size which are present at time $t$ have a finite
generation, i.e. result from finitely many dislocations of the initial particle. In this direction,
let us fix some arbitrarily small $\varepsilon>0$, and consider
the threshold
operator $\varphi_\varepsilon$ which consists of removing all the fragments with size less
than or equal to $\varepsilon$. Recall from Section 2.1
that $\varphi_\varepsilon(X)$ is a Markov chain, in particular the number of jumps accomplished
by this chain before time $t$ is finite a.s. This number obviously is an upperbound
for the generation of all fragments with size greater than $\varepsilon$. \QED

\subsection{Intrinsic martingale}

The purpose of this section is to introduce the so-called {\it intrinsic  martingale}
which is naturally induced by Malthusian hypotheses and the genealogical coding of
fragmentations, and  plays a crucial role in the asymptotic behavior of the latter.

\begin{theorem}\label{PSF3} The process
$${\cal M}_n:=\sum_{|u|=n}\xi_u^{p^*}\,,\qquad n\in\Z_+$$
is a martingale which is bounded in $L^1(\P)$,
and in particular, uniformly integrable. 
Moreover, the terminal value ${\cal M}_{\infty}$ is strictly positive a.s. whenever
$\nu(s_1=0)=0$ (i.e. a fragment may never disappear entirely after a dislocation).
\end{theorem}

In the sequel, $({\cal M}_n, n\in\Z_+)$  will be referred to as the {\rm intrinsic
martingale}. Observe that in the important
case when dislocations are conservative, in the sense that $\sum_{i=1}^{\infty}s_i= 1$,
$\nu(d{\bf s})$-a.s.,
then $p^*=1$ and ${\cal M}_n=1$ for all $n\in \Z_+$, and the statement is trivial.

 In general the distribution of the terminal value of the intrinsic martingale is not known
explicitly. However, it is straightforward from the branching property that there is the
identity in law 
\begin{equation}\label{eqFed}
{\cal M}_\infty\stackrel{\rm (d)}{=} \sum_{j=1}^{\infty} \xi_j^{p^*}{\cal M}^{(j)}_{\infty}
\end{equation}
where $\xi=(\xi_j, j\in\N)$ has the law $\nu$, and
${\cal M}^{(j)}_{\infty}$ are independent copies of ${\cal M}_\infty$, 
also independent of $\xi$. It is
known that under fairly general conditions, such equation characterizes the law
${\cal M}_\infty$ uniquely, see e.g. \cite{Liu2, R}. 
We also refer to Liu \cite{Liu1} for information of the tail behavior of the solution.

\section{Some applications}
In this section, we present some results on the behavior of self-similar fragmentation
chains which can be derived by the combination of standard techniques from the theories of
Markov chain in continuous time and multiplicative cascades.
\subsection{Some strong limit theorems $(\alpha\geq0)$}

The Malthusian parameter $p^*$ and the terminal value ${\cal M}_{\infty}$ of the intrinsic
martingale have a crucial role in the study of the asymptotic behavior of additive
functionals of the fragmentation, i.e. of the type
$$A(X(t))\,=\,\sum_{i=1}^{\infty}X_i^{p^*}(t)f(X_i(t),t)\,,\qquad t\geq0$$
for some measurable function $f: \R_+^2\to\R_+$. Indeed, 
let us denote by ${\cal M}(t)$ this quantity for $f\equiv 1$, so ${\cal M}(t)$ can be viewed as the analog of the intrinsic martingale
when the parameter is time instead of generation. It is easy to check that for $\alpha \geq 0$, ${\cal M}(t)$ converges a.s. and in $L^1(\P)$ to the terminal value ${\cal M}_{(t)}{\infty}$ of the intrinsic martingale. 
More generally, the branching property and a variation of the law of large number
enable to show that under appropriate
hypotheses on the function $f$, it holds that
\begin{equation}\label{eqmeta}
\lim_{t\to\infty} A(X(t))\,=\,c(f){\cal M}_{\infty}\,,\qquad
\hbox{in $L^1(\P)$}.
\end{equation}
Here, $c(f)$ is deterministic factor which can be determined by first
moment calculations. The terminal value of the intrinsic martingale ${\cal M}_{\infty}$
thus appears as a kind of universal random weight. We refer e.g. to \cite{Jag, Ner} and
references therein for many results in this vein for certain branching processes.
In order to avoid some technical discussion related to periodicity, we shall often make the
assumption that the dislocation law
$\nu$ is {\it non-geometric}, in the sense that there exists no real number 
$r\in]0,1]$ such that  $s_i\in\left\{r^n, n\in\Z_+\right\}\cup\{0\}$ for all $i\in\N$,
$\nu(d{\bf s})$-a.s.

Let us first consider the case of homogeneous fragmentation chains, i.e. with index of
self-similarity $\alpha=0$. In this situation, each particle has a standard exponential
lifetime, so informally when $n\in\N$ is large, particles at generation $n$
are alive at times close to $n$ (by the law of large numbers), and for the same
reason, when $t$ is large, the generation of particles alive at time $t$ is close to 
$[t]$. In this direction, one naturally expects that the strong limit theorems
for multiplicative cascades can be shifted without significant modifications to
homogeneous fragmentation chains. This is indeed the case as shown in the following result
which can be traced back (in a simpler setting) to Kolmogorov \cite{Ko}; see also
\cite{AK, Bi3, U} for closely related statements.

\begin{proposition}\label{TSF1}   Let $f:\R\to\R$ be a continuous
bounded function. In the homogenenous case $\alpha=0$, the following limits hold in
$L^1(\P)$:
$$\lim_{t\to\infty}\sum_{i=1}^{\infty}X_i^{p^*}(t) f(t^{-1}\ln X_i(t))
\,=\,{\cal M}_{\infty}f(-\kappa'(p^*))\,$$
and 
$$\lim_{t\to\infty}\sum_{i=1}^{\infty}X_i^{p^*}(t) f(t^{-1/2}(\ln X_i(t)+\kappa'(p^*)
t))\,=\,{\cal M}_{\infty}\E(f({\cal N}(0,-\kappa''(p^*)))\,,$$ 
where ${\cal N}(0,-\kappa''(p^*))$ denotes
a centered Gaussian variable with variance $-\kappa''(p^*)$. 
\end{proposition}

Roughly, the first part of Proposition \ref{TSF1} claims that in a homogenenous
fragmentation, `most' fragments decay exponentially fast with rate $\kappa'(p^*)$. The
second part is sharper and specifies the pathwise fluctuations. 

When the index of self-similarity $\alpha$ is positive, large fragments split faster
than small ones, so the rate of fragmentation decays as time passes and one may expect a
homogenisation phenomena. Theorem \ref{TSF2} below was first established by Filippov
\cite{Fi} in the special case when the dislocation measure $\nu$ is conservative, i.e.
when $\sum_{i=1}^{\infty}s_i=1$ for $\nu$-a.e. ${\bf s}$ (see also Brennan and Durrett
\cite{BD}). Recall that in this case, $p^*=1$ and ${\cal M}_n\equiv 1$.
The general case of non-conservative dislocation measures was recently proved by
Bertoin and Gnedin \cite{BG}.

\begin{theorem}\label{TSF2} 
Suppose that $\alpha>0$
and that the dislocation law $\nu$ is non-geometric. 
Then for every bounded continuous function $f:\R_+\to\R$
$$\lim_{t\to\infty}\sum_{i=1}^{\infty}X_i^{p^*}(t) f(t^{1/\alpha} X_i(t))
\,=\,{\cal M}_{\infty}\int_{0}^{\infty}f(y)\rho(dy)\,,\qquad \hbox{in }L^1(\P),$$
 where ${\cal M}_{\infty}$ is the
terminal value of the intrinsic martingale and 
$\rho$ is a deterministic probability measure.
More precisely, $\rho$ is determined by the moments
$$\int_{]0,\infty[}y^{\alpha k}\rho(dy)\,=\,{(k-1)!\over \alpha \kappa'(p^*)\,
\kappa(p^*+\alpha)\cdots\kappa(p^*+(k-1)\alpha)}\quad
\hbox{for $k\in\N$},
$$
(with the usual convention that the right-hand side above equals $1/(\alpha \kappa'(p^*))$
for
$k=1$).
\end{theorem}

Comparison with the homogeneous case of Proposition \ref{TSF1} shows that the size of
a typical fragment now decays as a power function of time. It is also interesting to
observe that the limit is much more sensitive to the dislocation law than
in the homogeneous situation: the function $\kappa$ can be recovered from the limit
measure $\rho$, whereas in the homogeneous case, the exponential rate of decay
just depends on the derivative of $\kappa$ at the Malthusian parameter.

Let us just sketch an argument for the proof of Theorem \ref{TSF2} based on moment
calculations.  Kolmogorov's backwards equation combined with Proposition \ref{P1} yields,
in the special case of a power function, the equation
$${d\over dt}\E\left(\sum_{i=1}^{\infty}X_i^{p-\alpha}(t)\right)
=-\kappa(p-\alpha)\E\left(\sum_{i=1}^{\infty}X_i^p(t)\right)\,.$$
The solution is given in the form
$$\E\left(\sum_{i=1}^{\infty}X_i^{p}(t)\right)\,
=\,\sum_{n=0}^{\infty}{(-t)^n\over n!}\gamma(n,p)\,,$$
where $\gamma(0,p)=1$ and for $n\geq1$
$$\gamma(n,p)\,=\,\prod_{k=0}^{n-1}\kappa(p+\alpha k)\,.$$
Asymptotics can be obtained using techniques from complex analysis, and one gets for
$p=p^*+\alpha k$ with
$k\in\N$, that
$$
\E\left(\sum_{i=1}^{\infty}X_i^{p}(t)\right)\sim
{(k-1)!\over \alpha \kappa'(p^*)\,
\kappa(p^*+\alpha)\cdots\kappa(p^*+(k-1)\alpha)} t^k\,,\qquad t\to\infty.
$$
This shows that
$$\lim_{t\to\infty}\E\left(\sum_{i=1}^{\infty}X_i^{p^*}(t) f(t^{1/\alpha} X_i(t))\right)
\,=\,\int_{0}^{\infty}f(y)\rho(dy)\,$$
when $f(x)=x^{k\alpha}$, and then for any continuous function $f$ bounded by a power
function. Finally, $L^1$-convergence in Theorem \ref{TSF2} follows from the meta-limit
theorem (\ref{eqmeta}). We refer to \cite{BG} for details.

\subsection{Extinction and formation of dust for $\alpha<0$}
Intuitively,  when the index of self-similarity is negative, fragments with
small sizes are  subject to high splitting rates, and this makes
them vanish entirely; see e.g. \cite{FG, H1, Je} for some works formalizing this
intuition.  

\begin{proposition}\label{PSF9} Suppose
$$\int_{\s} {\tt Card}\left\{i: s_i=1\right\}\nu(d{\bf s})\, < \, 1\,.$$
 Then the following assertions
hold  with
probability one:

\noindent {\rm (i)} For $\alpha<0$, 
$X(t)=(0,\ldots)$ for all sufficiently large $t$.

\noindent {\rm (ii)} When $\kappa(-\alpha)>0$, for almost every $t>0$
$${\tt Card}\left\{j\in\N: X_j(t)>0\right\}\,<\,\infty\,.$$
\end{proposition}

We stress that in general, no matter what
the value of $\alpha$ is, there may exist random instants 
$t$ at which 
$${\tt Card}\left\{j\in\N: X_j(t)>0\right\}\,=\,\infty\,.$$
For instance in the case when the dislocation law fulfills
$$
\nu\left(x_j>0 \hbox{ for all }j\in\N\right)\,=\,1\,,
$$
then with probability one, there occur infinitely
many sudden dislocations in the fragmentation
chain $X$, each of which produces infinitely many terms.
This does not induce any contradiction with Proposition \ref{PSF9} (ii)
when $\kappa(-\alpha)>0$, because informally, as the index of self-similarity is
negative, we know that fragments with small size vanish quickly.

It would be interesting to have information on the distribution of the extinction time
$$\zeta:=\inf\left\{t\geq0: X(t)=(0,0,\ldots)\right\}\,,$$
however it does not seem possible to express this law in a closed form.
Nonetheless, we point that an application of the branching property at the first dislocation
yields the identity in distribution
\begin{equation}\label{eqFed'}
\zeta\stackrel{\rm (d)}{=} {\bf e} + \max_{j\in\N} \xi_j^{-\alpha}\zeta_j\,,
\end{equation}
where ${\bf e}$ is a standard exponential variable, $(\xi_j, j\in\N)$ is distributed according to
$\nu$, $(\zeta_j, j\in\N)$ is a sequence of independent copies of $\zeta$, and finally,
${\bf e}$, $(\xi_j, j\in\N)$ and $(\zeta_j, j\in\N)$ are independent. We refer to \cite{AldB}
for a survey of this type of equations in distribution.

It is interesting to study in further details the extinction phenomenon in the case when
the dislocation measure is conservative, i.e. when
\begin{equation}\label{eqSF22}
\nu\left(\sum_{i=1}^{\infty} s_i\neq1\right)\,=\,0\,.
\end{equation}
It is easy to deduce by iteration that for every $n\in\N$, the total mass of particles
at the $n$-th generation is conserved, i.e.
$$\sum_{|u|=n}\xi_u\,=\,1\,,\qquad \hbox{a.s.}$$
Turning our interest to the total mass of particles at time $t$, we introduce the quantity
$$D(t):=1-\sum_{n=1}^{\infty}X_i(t)\,,$$
which can be viewed as the total mass of {\it dust}, that is of infinitesimal particles at time $t$.
One could be tempted to believe that the assumption (\ref{eqSF22}) would entail $D\equiv0$; 
indeed it is easy to check that this holds when the index of self-similarity of the
fragmentation is nonnegative. However
Proposition \ref{PSF9} shows that for negative indices of self-similarity,
$D$ reaches $1$ at a finite time a.s.

\begin{proposition}\label{PSFdust} Suppose (\ref{eqSF22}) holds. The following assertions
hold with probability one:

\noindent {\rm (i)} $D$ is a continuous
increasing process 
 which reaches $1$ in finite time. 

\noindent {\rm (ii)} If $\#(t):={\rm Card}\left\{i: X_i(t)>0\right\}$ denotes the number of fragments with
positive mass at time $t$, then
$$\int_{0}^{\infty}{\bf 1}_{\{\#(t)<\infty\}}dD(t)\,=\,0\,.$$
\end{proposition}

This statement again reflects the fact that, loosely speaking, dislocations occur
faster and faster as time passes. Observe that it entails that almost-surely, there exists uncountably
many times at which there are infinitely many fragments with positive size, which may be rather surprising
(for instance in the case when dislocations are binary, i.e. produce exactly two fragments; see
also Proposition \ref{PSF9}).

Recall that when the dislocation law $\nu$ is conservative, $\kappa(p)>0$ for every
$p>1$, and Proposition
\ref{PSF9} shows that for $\alpha<-1$, at each fixed time $t$ there is only a finite number
of fragments with positive size in the system. Combining this observation
with Proposition \ref{PSFdust} entails that the random measure $dD(t)$ is singular with
respect to Lebesgue measure on $\R_+$. Haas \cite{H2} established the following sharper
result.

\begin{proposition}\label{PSFhaas} Suppose (\ref{eqSF22}) holds and that there exists
$k\in\N$ such that
$\nu(s_{k+1}>0)=0$ (i.e. each dislocation splits a fragment into at most $k$ pieces).
The following assertions hold with probability one:

\noindent {\rm (i)} If $\alpha\in]-1,0[$, then $D$ is absolutely continuous. 

\noindent {\rm (ii)} If $\alpha<-1$, then the random measure $dD(t)$ is supported by a set
with Hausdorff dimension $1/|\alpha|$.
\end{proposition}

\section{Fragmentation chains and general branching processes}

We now turn our attention to a slightly different aspect of fragmentation chains,
by pointing at a connection with general (i.e. Crump-Mode-Jagers) branching processes; see
\cite{Jag, Ner} and references therein. Specifically, we use the genealogical coding of
the fragmentation chain of Section 2, and for each node $u\in{\cal U}$, we think of
$\sigma_u:=-\ln \xi_u$ as the
 birth-time of an individual. Every individual which is born has an infinite lifetime, and
gives birth to children according to a random point process
$\eta_i$. More precisely, if $s_1\geq s_2 \geq \ldots$ denote the sequence
ranked in the decreasing order of the sizes of the fragments that result
from the first dislocation of a fragment with size
$s$ (recall that we assume that $s\geq s_1$, i.e. dislocations always
produce fragments of smaller sizes), the $n$th child of that individual
is born at time $-\ln s_n$. In other words,
each fragment with size $s$ is viewed as an individual $u$  born at time
$\sigma_u=-\ln s$ and which has offspring described by the point process
$$\eta_u([0,t])\,=\,{\rm Card}\left\{n\in\N: s_n\geq\e^{-t}\right\}\,.$$

The particle system obtained in this way is then a general
branching process  with reproduction intensity
$$\mu(t)\,=\,\E(\eta_{\emptyset}([0,t]))\,=\,
\int_{\s}\sum_{n=1}^{\infty}{\bf 1}_{\{s_n\geq
\e^{-t}\}}\nu(d{\bf s})\,,\qquad t\geq0.$$
This point of view is useful to investigate problems related to situations
when  fragments
with size less than a certain fixed parameter $\varepsilon>0$ are
instantaneously frozen. The fragmentation process then terminates when the system only
constists in particles with size less than $\varepsilon$. 
This setting arises for instance in the mining industry where fragmentation is needed
to reduce  rocks into sufficiently small particles. For this
purpose, rocks are broken in crushers and mills by a repetitive
mechanism. Particles are screened so that when they become smaller than
the diameter of the mesh of a thin grid, they are removed from the
process. 

Imagine the instantaneous dislocation of a block of size $s$ into a set of smaller blocks
of sizes $(s_1,s_2,...)$ requires an energy 
of the form $s^\beta \varphi(s_1/s,s_2/s,\ldots)$, where $\varphi:\s\to\R$ is a  cost
function and $\beta >0$ a fixed parameter. We are interested in the total
energy-cost of the process that finishes when all particles have size less than
$\varepsilon>0$:
$${\cal E}(\varepsilon)\,:=\,
\sum_{u\in {\cal U}} {\bf 1}_{\{\xi_u> \varepsilon\}}
\xi_u^{\beta} \varphi(\tilde \xi_{u\bullet})\,,$$
where $\tilde \xi_{u\bullet}=(\xi_{u,n}/\xi_u, n\in\N)$
is the sequence of ratios of the sizes of the fragments resulting from
the dislocation of the fragment labelled by the node $u$ and the size of
that fragment (see Proposition \ref{PSF2'}).

The following asymptotic result for the energy cost as $\varepsilon\to0$ has been derived
recently by Bertoin and Martinez \cite {BM} from the work of Nerman \cite{Ner} on strong limit theorems for general branching processes. 

\begin{theorem}\label{TSF6} Suppose  that the dislocation law is
non-geometric and fulfils
$\int_{\s}
|\varphi({\bf s})|\nu(d{\bf s})<\infty$.
Then $\lim_{\varepsilon\to0}{\cal E}(\varepsilon):={\cal E}(0+)$ in $L^1(\P)$ for $\beta > p^*$, 
whilst for $\beta < p^*$ it holds that
$$\lim_{\varepsilon\to0}\varepsilon^{p^*-\beta}{\cal E}(\varepsilon)
={{\cal M}_{\infty} \over (p^*-\beta)\kappa'(p^*)}\int_{\s}
\varphi({\bf s})\nu(d{\bf s})\qquad \hbox{in }L^1(\P)\,.$$
\end{theorem}

 We may also consider the distribution of the
small particles that can go across the grid.
Specifically, we would like to get information about
the random finite measure
on $]0,1[$
$$\varrho_{\varepsilon}(dx)\,:=\,\sum_{u\in{\cal U}, u\neq\emptyset}
{\bf 1}_{\{\xi_{u-}\geq \varepsilon, \xi_u<\varepsilon\}}\xi_u^{p^*}
\delta_{\xi_u/\varepsilon}(dx)\,,$$
which can be viewed as a weighted version of the empirical
measure of the particles taken at the instant when they become smaller
than $\varepsilon$ and then renormalized.
In the setting of the general branching process associated to the fragmentation
chain, we can
re-express $\varrho_{\varepsilon}$ in the form
$$\langle \varrho_{\varepsilon},f\rangle
=\sum_{u\in{\cal U}}\sum_{n=1}^{\infty}\e^{-p^* \sigma_{un}}f(\e^{t-\sigma_{un}}){\bf
1}_{\{\sigma_u< t \leq \sigma_{un}\}}\,, 
$$
where $t=-\ln \varepsilon$ and $f:]0,1[\to\R_+$ denotes a generic measurable function.

\begin{proposition}\label{PSF4} Suppose that the dislocation law $\nu$ is non-geometric.
As $\varepsilon\to0$,  $\varrho_{\varepsilon}$ converges in
probability to ${\cal M}_\infty \mu$, where $\varrho$ is a deterministic probability
measure on $[0,1]$ given by
$$\varrho(dx)\,=\,\left(\int_{\s}\sum_{i=1}^{\infty}
{\bf 1}_{\{s_i<
x\}}s_i^{p^*}\nu(d{\bf s})\right) \,{dx\over x \kappa'(p^*)}\,.$$
\end{proposition}

\section{Fragmentations with instantaneous dislocations}

In this  section, we turn our attention to the situation when dislocations
can occur instantaneously. Examples arise naturally e.g. in the study of Brownian motion
and Continuum Random Trees; see \cite{AP, Be1, HM, Mie1, Mie2}.
Then one cannot consider the first
dislocation of a particle, and the genealogical structure of the process is no longer
discrete, which impedes the representation as a multiplicative cascade.
Our purpose is two-fold; first we would like to 
characterize such fragmentation processes and their
structure, and second we aim at extending properties which have been established
for fragmentation chains to this more general setting.

We shall focus on the case when the process $X$ takes values in the space of {\it
mass-partitions}
$$\p_{\rm m}:=\left\{{\bf s}\in\s: \sum_{i=1}^{\infty}s_i\leq 1\right\}\,.$$
A mass-partition ${\bf s}\in\p_{\rm m}$ can be thought of as the ranked sequence
of the masses of the fragments of some object with unit total mass; the case
when $\sum_{i=1}^{\infty}s_i< 1$ corresponds to the existence of dust, i.e.
a part of the object with positive mass has been reduced to infinitesimal particles each
with zero mass. 

It might be natural to try to investigate fragmentations with instantaneous dislocations
as limit of fragmentation chains when the parameter of the exponential lifetime
of a particle with unit mass tends to $\infty$ (so the lifetime of that particle tends to
$0$). However this approach is far from being easy, and we shall rather follow
a different route initiated by Kingman \cite{King} and then further developed by Pitman
\cite{Pit} in the framework of coalescents (see also Schweinsberg \cite{Sch}).

\subsection{Exchangeable partitions and interval representations}
Kingman's idea is to encode mass-partitions by random partitions of $\N$.
Specifically, consider some space $E$ endowed with a probability measure $\varrho$, and 
disjoint measurable subsets $E_1,\ldots$ with positive $\varrho$-measure. Set
$E_0:=E\backslash\left(\cup_{k\geq1}E_k\right)$, so $(E_0,E_1,\ldots)$ is a partition of
$E$. The sets $E_1,\ldots$ are viewed as fragments of $E$ and $E_0$ as the set of dust.

Then let $U_1, \ldots$ be a sequence of i.i.d. variables in $E$ with law $\varrho$,
and define a random partition $\pi$ of $\N$ as follows.
If $U_i\in E_0$, then $\{i\}$ is a singleton of the partition $\pi$,
and the other blocks of $\pi$ are of the type
$\pi_k=\left\{\ell\in \N: U_{\ell}\in E_k\right\}$ for $k\geq1$.
The strong law of large numbers enables us to recover the sequence of the masses
of the partition of $E$ as the asymptotic frequences of the blocks of $\pi$; more precisely
$$\varrho(E_k)\,=\,\lim_{n\to\infty}{1\over n}{\rm Card}\left\{\ell\in \pi_k: \ell\leq
n\right\}.$$

We stress that only a special class of partitions of $\N$ can arise in this framework, as
blocks of such partitions always have asymptotic frequences. More importantly, the random
partitions resulting from Kingman's construction are {\it exchangeable}, that is their
distribution is invariant by the action of permutations. Indeed the action of a
permutation  $\sigma$ amounts to a permutation of the indices
in the sequence of i.i.d. variables, and since $U_{\sigma(1)}, U_{\sigma(2)}, \ldots$
has the same law as $U_1, U_2, \ldots$, we see that the image of $\pi$
by the action of $\sigma$ has the same law as $\pi$.

Kingman \cite{King} (see also Aldous \cite{Aldsf} for a simpler proof) has shown
that for any random exchangeable partition $\gamma$ of $\N$, the blocks of $\gamma$ possess
asymptotic frequences, and more precisely, $\gamma$ has the same 
distribution as some partition $\pi$ constructed as above for a certain random
partition of a space $E$. Specifically, one can take $E=]0,1[$ endowed with the
Lebesgue measure, let
$E_1, \ldots$ be the interval components of some random open set $G\subseteq ]0,1[$
and $E_0=]0,1[\backslash G$.
Then $U_1,\ldots$ is a sequence of i.i.d. uniform variables which is independent of the
random open set $G$.

Now consider a fragmentation process $X$ in the sense of Definition \ref{D1}, and assume
that $X$ takes values in the space of mass-partitions. For our purpose, it is
convenient to  think of $X$ in terms of a fragmentation of the unit interval, in the sense
that there is a Markov process $(G(t), t\geq0)$ with values in the space of open sets in
$]0,1[$, such that $G(t)\subseteq G(s)$ when $s\leq t$ and for each $t\geq0$, $X(t)$ is
the ranked sequence of the interval components of $G(t)$. The fact that such a
representation exists is explained in
\cite{Ber1, Be3}.

Next, let $U_1,\ldots$ be a sequence of i.i.d. uniform variables which is independent
of $(G(t), t\geq0)$, and for each $t\geq0$, write $\Pi(t)$ for the random exchangeable
partition of $\N$ such that two distinct indices $i,j$ belong to the same block of
$\Pi(t)$ if and only if $U_i$ and $U_j$ belong to the same interval component of $G(t)$.
We make the key observation that for every $t\geq0$, given an interval component, 
say $I$ of $G(t)$, if $B=\{i\in \N: U_i\in I\}$ denotes the block of $\Pi(t)$
corresponding to $I$, then conditionally on $B$, $(U_i, i\in B)$ is a sequence of i.i.d.
uniform variables in $I$.
Essentially, this observation implies that the partition-valued process
$\Pi=(\Pi(t), t\geq0)$ is Markovian.
Thus to each fragmentation process $X$ we can associate a Markov process $\Pi$
with values in the space of partitions of $\N$ such that $\Pi(t)$ gets finer as $t$
increases. More precisely, the law of $\Pi$ is invariant by
the action of permutations and the branching property of $X$ is transferred to $\Pi$. Thus
$\Pi$ can be thought of as a fragmentation process with values in the space of partitions
of
$\N$.

Conversally, given an arbitrary fragmentation process $\Pi$ with values in the
space of partitions of $\N$ and which is invariant by the action of permutations, Kingman's
theorem enables us to recover a fragmentation process $X$ with values in $\p_{\rm m}$ by
considering the asymptotic frequences of the blocs of $\Pi(t)$.
We refer to \cite{Ber1} for a precise argument.

\subsection{The structure of homogeneous fragmentation processes}
In this section, we suppose that the index of self-similarity is $\alpha=0$.
In this situation, the fragmentation process
$\Pi$ with values in the space of partitions of $\N$ associated to $X$ as above enjoys a
crucial additional property : for every $n\in\N$, the restricted process $\Pi_{\mid [n]}$
is Markovian, where for a partition $\pi$ of $\N$, $\pi_{\mid [n]}$ denotes the
restriction of $\pi$ to $[n]:=\{1,\ldots, n\}$. Informally, this follows
from exchangeability and the fact that the dislocation rates of blocks of 
$\Pi$ do not depend on their sizes. Since the space of partitions of $[n]$ is finite,
 $\Pi_{\mid [n]}$ is a continuous time Markov chain, whose evolution is thus specified by
its jump rates. Combining this with Kingman's theory of exchangeable random partitions
enables us to reveal the structure of homogeneous fragmentation processes;
see \cite{Ber1,Be2} for details.
 Roughly, one gets that homogeneous fragmentations result from the combination
of two different phenomena: a continuous erosion and sudden dislocations. The erosion is a
deterministic mechanism, analogous to the drift for subordinators, whereas the
dislocations occur randomly according to some Poisson random measure, and can be viewed as
the jump-component of the fragmentation.

More precisely, let us first focus on dislocations. We call {\it dislocation measure} a
measure
$\nu$ on 
$\p_{\rm m}$, which gives no mass to $(1, 0,\ldots)$ and fulfills 
\begin{equation}\label{clev}
\int_{\p_{\rm m}}\left(1-x_1\right)\nu(dx)\,<\, \infty\,.
\end{equation}
Then, consider  a  Poisson
random measure on $\p_{\rm m}\times \N\times \R_+$,
$$\sum_{i=1}^{\infty}\delta_{(\Delta(i), k(i),t(i))}\,,$$
with intensity $\nu\otimes \#\otimes dt$, where $\#$ denotes the counting measure
on $\N$. One can construct a pure jump process $\left(Y(t), t\geq0\right)$
in $\p_{\rm m}$ which
jumps only at times $t(i)$ at which there is an atom $(\Delta(i),k(i),t(i))$ of the
Poisson measure. More precisely, the jump (i.e. the dislocation) induced by such an atom
can  be described as follows. 

The mass-partition $Y(t(i))$ at time $t(i)$ is obtained from that immediately
before $t(i)$, i.e. $Y(t(i)-)$, by replacing  its 
$k(i)$-th term, viz. $Y_{k(i)}(t(i)-)$, by the sequence
$Y_{k(i)}(t(i)-)\Delta(i)$, and ranking all the terms in the decreasing order. 
For instance, if 
$$Y(t(i)-)\,=\,\left({2\over 3},{1\over 4},{1\over 12},0,\ldots\right)\quad ,\quad
k(i)=2\quad 
\hbox{ and }\quad 
\Delta(i)\,=\,\left({3\over 4},{1\over 4},0,\ldots\right)\,$$ then we look at the 2-nd
largest term in the sequence
$Y(t(i)-)$, which is ${1\over 4}$, and split it according to $\Delta(i)$. This produces
two fragments of size ${3\over 16}$ and ${1\over 16}$, and thus
$$Y(t(i))\,=\,\left({2\over 3},{3\over 16},{1\over 12},{1\over 16},0,\ldots\right)\,.$$

Next, call {\it erosion coefficient} an arbitrary real number $c\geq0$, and set
$X(t)=\e^{-ct} Y(t)$.  
So $X$ is obtained from $Y$ by letting the fragments of the latter be eroded at constant
rate $c$. The process $(X(t), t\geq0)$ is again a homogeneous fragmentation,
and conversely any homogeneous fragmentation process can be constructed like this. In
conclusion, the distribution of a homogeneous fragmentation is entirely specified by its
erosion rate and its dislocation measure.

In the special case when $c=0$ and $\nu$ is a finite measure (and, in particular, a
probability), it is easy to check that $X$ is then a homogeneous fragmentation chain with
dislocation law
$\nu$ in the sense of Section 2. 
Note however that the condition (\ref{clev}) allows 
$\nu$ to be infinite, which corresponds to the situation when particles dislocate
instantaneously. Informally, mass-partitions ${\bf s}$ for which $1-s_1$ is small should
be thought of as small, in the sense that a small mass-partition produces one large
fragment and all the remaining ones are small (in particular, the mass-partition
$(1,0,\ldots)$ has to be viewed as a neutral element for dislocations). So, roughly
speaking,  condition (\ref{clev}) allows infinitely many small dislocations, but guaranties
that the accumulation of these small dislocations does not reduce instantaneously
the initial mass into dust. This bears obvious similarities with subordinators,
which are constructed by the It\^o-L\'evy decomposition from atoms of a certain Poisson
random measure on $\R_+$, see Chapter 1 in \cite{Besf}. The intensity of this Poisson
random measure is given by the so-called L\'evy measure $\Lambda$ of the subordinator, and
the integrability condition $\int_{\R_+}(1\wedge x)\Lambda(dx)<\infty$ for a measure on
$\R_+$ to be the L\'evy measure of some subordinator is the necessary and sufficient
condition for the summability of the atoms.

The Poissonian structure of homogeneous fragmentation is a fundamental tool which
enables to circumvent difficulties related to the absence of a discrete genealogical
structure. However, although the law of a homogeneous fragmentation process is
characterized by its erosion rate and dislocation measure, in general we do not know how to
describe explicitly e.g. the distribution of the process at a fixed time. 
The next crucial tool for the study lies in the fact that 
partial but most useful information can be derived from a so-called size-biased
sampling, and it turns out that the law of the latter is simple to formulate.

Recall that we may represent $X$ in terms of fragmentation of the unit interval.
So we consider a Markovian family 
$(G(t), t\geq0)$ of nested open subsets of the unit interval, 
in particular for every $t\geq0$,  $X(t)$ is the ranked sequence of the lengths of the
interval components of $G(t)$. Now suppose that $U$ is a uniform random variable
on $]0,1[$, which is independent of $(G(t), t\geq0)$, and for every $t\geq0$, denote
by $\chi(t)$ the length of the interval component of $G(t)$ which contains $U$
(with the convention that $\chi(t)=0$ if $U\not\in G(t)$).
In other words, the process $\chi=(\chi(t), t\geq0)$ gives the size of the fragment
containing a point which has been tagged independently of the fragmentation process;
one refers to $\chi$ as the process of the tagged fragment.
Note that $\chi(t)$ is a size-biased pick from the sequence $X(t)=(X_1(t),\ldots)$, 
i.e. there is the identity in law
$$\chi(t)\mbox{$ \ \stackrel{\cal L}{=}$ }X_K(t)\,,$$
where $K$ is an integer valued variable
whose conditional distribution  given $X(t)$ is 
$$\P(K=k\mid X(t))\,=\,X_k(t)\,,\quad k=1,\ldots\,.$$
It turns out that the process of the tagged fragment
is closely related to a subordinator (i.e. an increasing process with independent and
stationary increments; see \cite{Besf} for background) that can be characterized explicitly
in terms of the dislocation measure and the erosion coefficient.
\begin{theorem}\label{THF1} The process
$$\sigma(t)\,=\,-\ln \chi(t)\,,\qquad t\geq0$$
  is a subordinator, whose law is determined by its Laplace transform of its
one-dimensional distributions. We have
$$\E\left(\chi(t)^q\right)\,=\,\E(\exp(-q\sigma(t)))\,=\,\exp\left(-t\kappa(q+1)\right)\,,\qquad
t,q >0\,,
$$
where the function $\kappa$ is given in terms of the erosion rate $c$ and the dislocation
measure
$\nu$ by the identity
\begin{equation}\label{equ1}
\kappa(q)\,:=\,cq + \int_{\p_{\rm
m}}\left(1-\sum_{i=1}^{\infty}s_i^q\right)\nu(d{\bf s})\,,\qquad q\geq 1\,.
\end{equation}

\end{theorem}

Many features (like asymptotic behavior) of fragmentation process can be read on properties
of tagged fragments, the combination of Theorem \ref{THF1} and the theory of subordinators
provides the key to many results on homogeneous fragmentations.

\subsection{Additive martingales and applications}

The independence and stationarity of the increments of subordinators entail that
for every $q\geq0$, the process
$$\exp(-q\sigma(t)+\kappa(q+1)t)\,=\,\exp(t\kappa(q+1))\chi^q(t)\,,\qquad t\geq0$$
is a martingale. As the tagged fragment $\chi(t)$ is picked at random from the
mass-partition $X(t)$ by size-biased sampling, it follows (take $q=p+1$) that
$$M(p,t):=\exp(t\kappa(p))\sum_{i=1}^{\infty}X_i^p(t)\,,\qquad t\geq0$$
is a nonnegative martingale for every $p>\underline p$
(this can also be checked directly from the branching and scaling properties). 

One refers to $M(p,t)$ as an {\it additive martingale}.
Plainly $M(p,t)$  converges a.s. as $t\to\infty$, and in order to investigate the
asymptotic behavior of homogeneous fragmentation processes, it is important to know if the
limit of this martingale is strictly positive or zero. 
We shall investigate this question in the special case when
the dislocation measure is conservative and there is no erosion, i.e.
we assume from now on that
\begin{equation}\label{eqH}
c=0\hbox{ and }\sum_{i=1}^{\infty}s_i=1\quad \nu(d{\bf
s})\hbox{-a.e.}
\end{equation} Observe that the Malthusian hypotheses are then automatically fulfilled 
and the Malthusian parameter is $p^*=1$.
A first step in the analysis is the following elementary lemma.

\begin{lemma}\label{LSF5} Assume (\ref{eqH}). The  function
$p\to\kappa(p)/p$ reaches its maximum at a unique location 
$\bar p>1$, which is the unique solution to the equation
$$p\kappa'(p)=\kappa(p)\,.$$
More precisely, the function $p\to\kappa(p)/p$ increases on $]\underline p,\bar p[$
and decreases on $]\bar p, \infty[$, and the value of its maximum is
$\kappa'(\bar p)=\kappa(\bar p)/\bar p$.

\end{lemma}
\proof 
We first point out that the function $\kappa$ is concave and increasing.
It follows that
\begin{equation}\label{eqSF19}
\hbox{the function
 $p\to p\kappa'(p)-\kappa(p)$ decreases on $]\underline p,\infty[$.}
\end{equation}
Indeed, this function has derivative $p\kappa''(p)$, which is negative
since $\kappa$ is concave. Recall from (\ref{equ1}) and (\ref{eqH}) that $\kappa(1)=0$; on
the other hand, it is obvious that $\lim_{q\to\infty}\kappa(q)/q=0$,  hence the function
$p\to\kappa(p)/p$ has the same limit at $1$ and at
$\infty$, so it  reaches its overall maximum  at a unique point
$\overline p>1$. In particular, we deduce from (\ref{eqSF19}) that the
derivative of $p\to\kappa(p)/p$ is positive on
$]\underline p,\overline p[$ and negative on
$]\overline p,\infty[$. 
Finally, the derivative must be zero at $\overline p$, which entails that
the overall maximum is given by $\kappa'(\overline p)=\kappa(\overline p)
/\overline p$. \QED

We may now state the main result of this section which can be proved using
the Poissonian structure of homogeneous fragmentations and stochastic calculus, see
\cite{Be4} for details.

\begin{theorem}\label{TSF3} Assume (\ref{eqH}). For every 
$p\in]\underline p,\bar p[$, the martingale $M(p,\cdot)$ is bounded in $L^1(\P)$ and
its terminal value is strictly positive.
\end{theorem}

Just as Theorem \ref{PSF3} for the intrinsic martingale, Theorem \ref{TSF3} has crucial
role in the study of the asymptotic behavior of homogeneous fragmentation. First,
we specify the rate of decay of the largest fragment, refering to Biggins \cite{Bi2}
for a similar result in the framework of branching random walks.

\begin{corollary}\label{CSF2}  Assume (\ref{eqH}). It holds with probability one that
$$\lim_{t\to\infty}{1\over t}\ln X_1(t)\,=\,-\kappa'(\bar p)
\,=\,-{\kappa(\bar p)\over \bar p}\,.$$
\end{corollary}

\proof For every $p>\underline p$, we have
$$\exp(t\kappa(p))X_1^p(t)\leq \exp(t\kappa(p))\sum_{i=1}^{\infty}
X_i^p(t)$$
and the right-hand side remains bounded as $t$ tends to infinity.
Hence
$$\limsup_{t\to\infty}{1\over t}\ln X_1(t)\,\leq\,-{\kappa( p)\over 
p}\,,$$ and optimizing over $p$ yields
$$\limsup_{t\to\infty}{1\over t}\ln X_1(t)\,\leq\,
-{\kappa( \bar p)\over \bar p}\,.$$

On the other hand, for every  $p\in]\underline p,\bar p[$ and $\varepsilon>0$
sufficiently small, we have the lower bound
$$\exp(t\kappa( p))\sum_{i=1}^{\infty}
X_i^p(t)\leq X_1^{\varepsilon}(t)\exp(t\kappa( p))\sum_{i=1}^{\infty}
X_i^{ p-\varepsilon}(t)\,.$$
We know that both limits
$$\lim_{t\to\infty}\exp(t\kappa( p))\sum_{i=1}^{\infty}
X_i^{ p}(t) \quad \hbox{and} \quad 
\lim_{t\to\infty}\exp(t\kappa( p-\varepsilon))\sum_{i=1}^{\infty}
X_i^{ p-\varepsilon}(t) $$
are finite and strictly positive a.s., and we deduce that
$$\liminf_{t\to\infty}{1\over t}\ln X_1(t)\,\geq\,
-{\kappa( p)-\kappa( p-\varepsilon)\over \varepsilon}\,.$$
We take the limit of the right-hand side as $\varepsilon\to0+$ and then
as $p$ tends to $\bar p$ to conclude that
$$\liminf_{t\to\infty}{1\over t}\ln X_1(t)\,\geq\,
-\kappa'( \bar p)\,.$$
Now, this quantity coincides with $-\kappa(\bar p)/\bar p$, as we know
from Lemma \ref{LSF5}.
\QED

It is interesting to compare Corollary \ref{CSF2} with Proposition \ref{TSF2}(i), which
claims that the size of most fragments decays exponentially fast with rate
$\kappa'(p^*)=\kappa'(1)$. The size of largest fragment thus also decays exponentially
fast, but with a slower rate $\kappa'(\bar p)<\kappa'(1)$. We refer to Berestycki
\cite{Ber2} for the multi-fractal analysis of the exponential rates of decay of fragments
in homogeneous fragmentations.

Alternatively, one can also establish Theorem \ref{TSF3} by
discretization, using the following connection with branching random walks and classical
results on the latter; see \cite{BR}. If we consider the
 point measure $Z^{(t)}$ with atoms at the logarithms of the 
fragments
\begin{equation}\label{defZ}
Z^{(t)}\,:=\,\sum_{i=1}^{\infty}\delta_{\ln
X_i(t)}\,,\qquad t\geq0\,,
\end{equation}
we can think of the discrete skeleton $(Z^{(n)}, n\in\N)$ as a non-interacting particle
system. Specifically, particles evolve independently one of the other, and at each step,
each particle, say
$y$, is replaced birth to a random cloud of particles $y+{\cal Z}$,
where the law of ${\cal Z}$ is that of $(\ln X_i(1), i\in\N)$. This means that $(Z^{(n)},
n\geq0)$ is a branching random walk in the sense of \cite{AK}, \cite{Bi1}, ... The latter
have been throughoutly studied in the literature, and many of their properties can be
translated to homogeneous fragmentation processes.
For instance, we can derive precise information on almost sure large
deviations for the empirical measure; using a genuine result of Biggins \cite{Bi3}
for branching random walks.

\begin{corollary}\label{c1}
 Assume (\ref{eqH}) and that 
the dislocation measure $\nu$ is non-geometric, and for 
 $p\in]\underline p, \bar p[$, let $M(p,\infty)$ denote the terminal
value of the uniformly integrable martingale $M(p,\cdot)$.
 If $f:\R\to\R$ is a function with compact
support which is directly Riemann integrable,
then
$$\lim_{t\to\infty}  \sqrt t \, 
\e^{-t(p\kappa'(p)-\kappa(p))}\int_{\R}^{}f(t
\kappa'(p)+y)Z^{(t)}(dy)\,=\, {M(p,\infty)\over \sqrt{2\pi
|\kappa''(p)|}}\int_{-\infty}^{\infty}f(y)\e^{-py}dy\,.
$$
uniformly for $p$ in compact subsets of $]\underline p, \bar p[$, almost
surely.
\end{corollary}

In a different direction, one can also use the discretization techniques to
estimate the probability of presence of abnormally large fragments as
time goes to infinity; see \cite{BR} for details.

\begin{corollary} Assume (\ref{eqH}) and that the dislocation measure $\nu$ is
non-geometric, fix two real numbers
$a<b$ and take  $p>\bar p$. Then as $t\to\infty$
$$ \P \left(\exists i\in \N: X_i(t)\in [ a\e^{-t\kappa'(p)}  , 
b\e^{-t\kappa'(p)}]\right)\,
\sim\, K t^{-1/2}  \e^{t(p\kappa'(p)-\kappa(p))} \,,$$
where $K$ is some positive and finite constant depending on
$a,b$ and the characteristics of the fragmentation.
\end{corollary}

\subsection{Changing the index of self-similarity}
So far, we have only been able to study fragmentations with instantaneous dislocations in
the homogeneous case, i.e. when the index of self-similarity is $\alpha=0$. In this 
section, we present a simple transformation that changes a homogeneous
fragmentation $X$ into a self-similar one with an arbritrary index of self-similarity,
$X^{(\alpha)}$, completing the construction of general self-similar fragmentation
processes.  We refer to \cite{Be3} for details. 

In this direction, it is convenient to start from
an interval-representation $(G(t), t\geq0)$ of some homogeneous fragmentation $X$ as in
section 6.1.  For every
$y\in]0,1[$, let
$I_y(t)$ denote the interval component of $G(t)$ that contains $y$ if $y\in G(t)$, and
$I_y(t)=\emptyset$ otherwise. We write $|I|$ for the length of an interval $I\subseteq ]0,1[$, and
for every $y\in]0,1[$ we consider the time-substitution
$$T^{(\alpha)}_y(t)\,:=\,\inf\left\{u\geq0: \int_{0}^{u}|I_y(v)|^{-\alpha}
dv>t\right\}\,.$$
Because the open sets $G(t)$ are nested, we see that for every $y,z\in]0,1[$, the intervals 
$I_y(T^{(\alpha)}_y(t))$ and $I_z(T^{(\alpha)}_z(t))$ are either identical or disjoint, so the family 
$\left\{I_y(T^{(\alpha)}_y(t)), 0< y <1\right\}$ can  be viewed as the interval components of an open
set $G^{(\alpha)}(t)$. It is straightforward that the family $(G^{(\alpha)}(t),t\geq0)$ is nested. More precisely,
if we write $X^{(\alpha)}(t)$ for the ordered sequence of the lengths of the interval components of
$G^{(\alpha)}(t)$, then $(X^{(\alpha)}(t), t\geq0)$ 
is a self-similar fragmentation with index $\alpha$. 

Any self-similar fragmentation $X^{(\alpha)}$ can be
constructed from some homogeneous one $X$ as above, and this construction can be inverted.
 In particular, the distribution of $X^{(\alpha)}$ is entirely
determined by the index of self-similarity
$\alpha$, and the erosion coefficient $c\geq0$ and the dislocation measure $\nu$ of the
homogeneous fragmentation $X$. 

The key tool which is needed to extend the results of Sections 3-5 to self-similar
fragmentation with instantaneous dislocations is provided by the stochastic structure
of the process of the tagged fragment, $(\chi^{(\alpha)}(t), t\geq0)$.
Recall that in the homogeneous case, $(\chi(t), t\geq0)$  can be described as the
exponential of a subordinator. The construction above of a self-similar fragmentation from
a homogeneous one by time substitution enables us to derive the law of the tagged fragment
$\chi^{(\alpha)}(t)$ for a self-similar fragmentation.

Specifically, let $\sigma=\left(\sigma(t), t\geq0\right)$ be a subordinator
with Laplace exponent $\kappa(1+\cdot)$ given by {\rm (\ref{equ1})}. Introduce the
time-change
$$\tau(t)\,=\,\inf\left\{u: \int_{0}^{u}\exp(\alpha\sigma(r))dr>t\right\}\,,\qquad
t\geq0\,,$$ and set $\zeta_t=\exp(-\sigma(\tau(t)))$ (with the convention that $\zeta_t=0$
if
$\tau(t)=\infty$). Then the processes $(\zeta_t, t\geq0)$ and $\left(\chi^{(\alpha)}(t),
t\geq0\right)$ have the same law. In particular, this shows that the process of the tagged
fragment  is a decreasing self-similar Markov process as introduced by
Lamperti \cite{La}, see the survey \cite{BY} and the references therein.

Results stated for self-similar fragmentation chains in Sections 4-5 can be extended
verbatim to self-similar fragmentation processes; see \cite{Be4, BG, BM, H2}.

\section{Duality with certain coalescent processes}

Coagulation processes are used as models to describe the evolution of particle systems in which pairs (or, more generally, families) of particles merge as time passes.
At this level of generality, a simple time-reversal provides an obvious connection with fragmentation.

However in practice, just as for fragmentation, one has to make some restrictive assumptions on the dynamics of coagulation in order to deal with
processes that can be studied mathematically. First, one 
assumes that particles are determined by their masses (i.e. a positive real number such that 
the total mass is a preserved quantity when merging occurs). So there is no geometry involved in the system; physicists call such models {\it mean-field}. Second, one assumes that the evolution  is Markovian, and third, that the rate at which
a family of particles merges only depends on the particles in this family, and not on the other particles in the system. We shall call {\it stochastic coalescent} a coagulation process that fulfils these requirements, and refer to the survey by Aldous
\cite{Ald} and the references therein for much on this notion.

The first two requirements (mean-field and Markov properties) are clearly compatible with time-reversal. However, even though the third one bears some vague resemblance with the branching property, there is in general no reason why it should yield the latter after time-reversal. Despite the absence of a general result on duality by time-reversal between fragmentation and coalescent processes, there are nontheless several important examples for which duality holds. It would be very interesting to establish general critera for duality;
we leave this question open as a challenge, and now conclude this survey by discussing some examples.

The first example is the simplest, it can be constructed from a sequence $U_{1},\ldots$
of independent uniform variables on $[0,1]$ and an independent Poisson process $N=(N_{t}, t\geq 0)$. Specifically, for each time $t\geq 0$, write $F(t)$ for the ranked sequence of the 
lenghts of the interval components of the random open set $]0,1[\backslash
\{U_{i}: 1\leq i \leq N_{t}\}$. It is easy to see that $F=(F(t), t\geq 0)$
is a self-similar fragmentation chain with index $\alpha=1$.
More precisely, its dislocation law is that of the random mass partition
$(1-U/2, U/2, 0,\ldots)$ where $U$ is uniformly distributed on $[0,1]$.
Then it can be checked  (see \cite{BeRSA} for details) that the exponential time reversal
$C(t):=F(\e^{-t})$ transforms the fragmentation process $F$ into a coalescent process $C$.
Specifically, when there are $n\geq 2$ particles in the coalescent $C$, the first coagulation
occurs after an exponential time with parameter $n$, and the pair of particles involved
is uniformly distributed amongst the $n(n-1)/2$ possible pairs. One reckognizes a variation
of the celebrated coalescent of Kingman \cite{King}  (which is used to describe the genealogy of large populations and  has many applications 
in Biology).  More generally, a similar duality holds for self-similar fragmentation chains
with index $\alpha=1$ and dislocation law given in terms of certain $n$-dimensional Dirichlet distributions. These fragmentation and coalescent processes appear in the genealogy of Yule processes; see \cite{BeGo}.

The second example concerns the {\it additive coalescent}, a coagulation process which arises e.g. as a model for the formation of drops of rain in clouds. Roughly speaking, in an additive coalescent, any pair of particles, say $(x,y)$, merges at rate $x+y$, independently of the other pairs in the system. Evans and Pitman \cite{EvPit} have observed that if $C^{+,n}=(C^{+,n}(t), t\geq 0)$ denotes the process started from the monodisperse initial condition which consists in $n$ particles each with mass $1/n$, then as $n\to \infty$, the process $(C^{+,n}(t+{1\over 2}\ln n), t\geq -{1\over 2}\ln n)$
converges in the sense of finite dimensional distributions to 
$C^{+}:=(C^{+}(t), t\in\R)$. The latter is known as
the {\it standard} additive coalescent. Again the exponential time-change, $F(t)=C^+(-\ln t)$, transforms the coalescent into a fragmentation process which is self-similar with index
$\alpha = 1/2$. We refer to \cite{AP, Be1} for more on this topic. In this direction,
we also mention that Miermont \cite{Mier0} has shown that the exponential time-change also transforms some (but not all) non-standard additive coalescent into fragmentation processes; however the latter are not self-similar in general.

The final example is due to Pitman \cite{Pit} who established a remarkable duality between
certain coagulation and fragmentation operators based on adequat Poisson-Dirichlet distributions. An important special case involves Ruelle cascades and the Bolthausen-Sznitman coalescent. The former have been introduced in \cite{Ruelle} as a tool for studying Derrida's Generalized Random Energy Model of spin glass; and the latter in \cite{BS} to describe the time-reversed dynamics in Ruelle's cascades.
The Bolthausen-Sznitman coalescent also appears naturally in the genealogy
of Neveu's branching process, see \cite{BLG}.
Pitman's duality result enables us to view Ruelle's cascades as a time-inhomogeneous
fragmentation process, see also Basdevant \cite{Bas} for a recent development in this direction.

\end{document}